\documentstyle{amsppt}
\NoRunningHeads

\def\r2{\Bbb R^2}
\def\z2{\Bbb Z^2}

\def\n{\Bbb N}
\def\q{\Bbb Q}
\def\p{\partial BX}
\def\ve{\varepsilon}

\def\eps{\epsilon}

\def\ff{\Bbb F}
\def\aa{\alpha}
\def\bb{\beta}
\def\qq{\Bbb Q}
\def\rr{\Bbb R}
\def\zz{\Bbb Z}

\topmatter
\author{Sergei Konyagin and Izabella {\L}aba}
\endauthor
\title
Distance sets of well-distributed planar sets for polygonal norms
\endtitle
\address Department oF Mechanics and Mathematics, Moscow State University, 
Moscow, 119992, Russia,
e-mail: konyagin\@ok.ru
\endaddress

\address Department of Mathematics, University of British Columbia,
Vancouver, B.C. V6T 1Z2, Canada, e-mail: ilaba\@math.ubc.ca
\endaddress
\abstract
Let $X$ be a 2-dimensional normed space, and let $BX$ be the unit ball in 
$X$.  We discuss the question of how large the set of extremal points of $BX$
may be if $X$ contains a well-distributed set whose distance set $\Delta$
satisfies the estimate $|\Delta \cap[0,N]|\leq CN^{3/2-\eps}$.  We also give 
a necessary and sufficient condition for the existence of a well-distributed 
set with $|\Delta\cap [0,N]|\leq CN$.
\endabstract
\endtopmatter
\document

\centerline{\S0. INTRODUCTION}
\bigskip
The classical Erd\H os Distance Problem asks for the smallest possible 
cardinality of
$$\Delta(A)=\Delta_{l_2^2}(A)=\left\{\|a-a'\|_{l_2^2}:a,a'\in A\right\}$$
if $A\subset\r2$ has cardinality $N<\infty$ and
$$\|x\|_{l_2^2}=\sqrt{x_1^2+x_2^2}$$
is the Euclidean distance between the points $a$ and $a'$.
Erd\H os conjectured that $|\Delta(A)|\gg N/\sqrt{\log N}$ for $N\ge2$.
(We write $U\ll V$, or $V\gg U$, if the functions $U,V$ satisfy the inequality
$|U|\le CV$, where $C$ is a constant which may depend on some specified 
parameters). The best known result to date in two dimensions is due to Katz and Tardos
who prove in \cite{KT04} that $|\Delta(A)|\gg N^{.864}$ improving an
earlier breakthrough by Solymosi and T\'oth \cite{ST01}. 

More generally, one can examine an arbitrary two-dimensional space $X$
with the unit ball
$$BX=\{x\in\r2:\,\|x\|_X\le1\}$$ 
and define the distance set
$$\Delta_{X}(A)=\left\{\|a-a'\|_X:a,a'\in A\right\}.$$
For example, let 
$$\|x\|_{l_\infty^2}=\max(|x_1|,|x_2|)$$
then for $N\ge1$, $A=\{m\in\z2:0\le m_1\le N^{1/2},0\le m_2\le N^{1/2}\}$
we have $|A|\gg N$, $|\Delta_{l_\infty^2}(A)|\ll N^{1/2}$. This simple example 
shows that the Erd\H os Distance Conjecture can not be directly extended for
arbitrary two-dimensional spaces. Erd\H os \cite{E46} (see also \cite{I01})
proved the estimate $|\Delta_{X}(A)|\gg N^{1/2}$ for any space $X$.

Also, for a positive integer $N$ we denote
$$\Delta_{X,N}(A)=\left\{\|a-a'\|_X\le N:a,a'\in A\right\}.$$

We say that a set $S\subset X$ is well-distributed if there is a constant 
$K$ such that every closed ball of radius $K$ in $X$ contains a point from 
$S$. In other words, for every point $x\in X$ there is a point $y\in S$ such
that $\|x-y\|_X\le K$. Sometimes it is said that $S$ is a $K$-net for $X$.  
Clearly, for any well-distributed set $S$ and $N\ge 2K$ we have
$$|\{x\in S:\|x\|_X\le N/2\}|\gg N^2\tag1$$ 
where the constant in $\gg$ depends only on $K$. Therefore, for any 
well-distributed set $S\in l_2^2$ we have, by \cite{T02}, 
$$|\Delta_{l_2^2,N}(S)|\gg N^{1.728},$$
and the Erd\H os Distance Conjecture implies for large $N$
$$|\Delta_{l_2^2,N}(S)|\gg N^2/\sqrt{\log N}.$$
On the other hand, for a well-distributed set $S=\z2\subset l_\infty^2$
we have
$$|\Delta_{l_\infty^2,N}(S)|=2N+1.$$

Iosevich and the second author \cite{I{\L}03} have recently established that a 
slow growth of
$|\Delta_{X,N}(S)|$ for a well-distributed set $S\subset X$ is possible only
in the case if $BX$ is a polygon with finitely or infinitely many sides. Let 
us discuss possible definitions of polygons with infinitely many sides. For a 
convex set $A\subset X$ by $Ext(A)$ we denote the set of extremal points of 
$A$. Namely, $x\in Ext(A)$ if and only if $x\in A$ and for any segment $[y,z]$
the conditions $x\in[y,z]\subset A$ imply $x=y$ or $x=z$. Clearly, $Ext(BX)$
is a closed subset of the unit circle 
$$\p=\{x\in X:\|x\|_X=1\}.$$ 
Also, it is easy to see that $Ext(BX)$ is finite if and only if $BX$ is
a polygon with finitely many sides, and it is natural to consider $BX$ 
as a polygon with infinitely many sides if $Ext(BX)$ is small. There are 
different ways to define smallness of $Ext(BX)$ and, thus, polygons 
with infinitely many sides:\newline
1) in category: $Ext(BX)$ is nowhere dense in $\p$;\newline
2) in measure: $Ext(BX)$ has a zero linear measure (or a small Hausdorff 
dimension);\newline
3) in cardinality: $Ext(BX)$ is at most countable.\newline
Clearly, 3) implies 2) and 2) implies 1).

It has been proved in \cite{I{\L}03} that the condition
$$\underline\lim_{N\to\infty}|\Delta_{X,N}(S)|N^{-3/2}=0\tag0.1$$
for a well-distributed set $S$ implies that $BX$ is a polygon in a 
category sense. Following \cite{I{\L}03}, we prove that, moreover, $BX$ is 
a polygon in a measure sense.
\proclaim{Theorem 1} Let $S$ be a well-distributed set.\newline
(i) Assume that (0.1) holds. Then the one-dimensional Hausdorff measure of 
\newline
$Ext(BX)$ is $0$;\newline    
(ii) If moreover 
$$|\Delta_{X,N}(S)|=O(N^{1+\alpha})\tag0.2$$
for some $\alpha\in(0,1/2)$ then the Hausdorff dimension of $Ext(BX)$ is
at most $2\alpha$.
\endproclaim

If $|\Delta_{X,N}(S)|$ has an extremally slow rate of growth for some
well-distributed set $S$, namely,
$$|\Delta_{X,N}(S)|=O(N)\tag0.3$$
then, as it has been proved in \cite{I{\L}03}, $BX$ is a polygon with finitely 
many sides. However, if we weaken (0.3) we cannot claim that $BX$ is
a polygon in a cardinality sense.
\proclaim{Theorem 2} Let $\psi(u)$ be a function $(0,\infty)\to(0,\infty)$
such that $\lim_{u\to\infty}\psi(u)=\infty$. Then there exists a 
$2$-dimensional space $X$ and a well-distributed set $S\subset X$ such that
$$|\Delta_{X,N}(S)|=o(N\psi(N))\quad(N\to\infty)\tag0.4$$
but $Ext(BX)$ is a perfect set (and therefore is uncountable).   
\endproclaim

Also, we find a necessary and sufficient condition for a space $X$
to make (0.3) possible for some well-distributed set $S\subset X$.
Take two non-collinear vectors $e_1,e_2$ in $X$. They determine coordinates
for any $x\in X$, namely, $x=x_1e_1+x_2e_2$. Then, for any non-degenerate
segment $I\subset X$, we can define its slope $Sl(I)$: if the line
containing $I$ is given by an equation $u_1x_1+u_2x_2+u_0=0$, then
we set $Sl(I)=-u_1/u_2$. We write $Sl(I)=\infty$ if $u_2=0$; it will be
convenient for us to consider $\infty$ as an algebraic number. 
\proclaim{Theorem 3} The following conditions on $X$ are equivalent:\newline
(i) $BX$ is a polygon with finitely many sides, and there is a coordinate 
system in $X$ such that the slopes of all sides of $BX$ are algebraic;\newline
(ii) there is a well-distributed set $S\subset X$ such that (0.3) holds.
\endproclaim
\proclaim{Corollary 1} If a norm $\|\cdot\|_X$ on $\r2$ is so that 
$BX$ is a polygon with finitely many sides and all angles between its sides
are rational multiples of $\pi$ then there is a well-distributed set 
$S\subset X$ such that (0.3) holds.
\endproclaim
\proclaim{Corollary 2} If a norm $\|\cdot\|_X$ on $\r2$ is defined by a 
regular polygon $BX$ then there is a well-distributed set 
$S\subset X$ such that (0.3) holds.
\endproclaim

%We say that a set $S\subset X$ is separated if 
%$$\inf\Sb a,a'\in S\\ a\neq a'\endSb |a-a'|>0.$$
%By the way we prove using arguments from \cite{IL03} the following.
%\proclaim{Proposition 1} Every well-distributed set satisfying (0.3)
%is separated.
%\endproclaim

The Falconer conjecture (for the plane) says that if the Hausdorff
dimension of a compact $A\subset\r2$ is greater than $1$ then $\Delta(A)$
has positive Lebesgue measure. The best known result is due to Wolff who proved
in \cite{W99} that the distance set has positive Lebesgue measure if the  
Hausdorff dimension of $A$ is greater than $4/3$. One can ask a similar 
question for an arbitrary two-dimensional normed space $X$. It turns out that
this question is related to distance sets for well-distributed and separated 
sets. By Theorem 4 from \cite{I{\L}04}, Theorem 3 and Proposition 1 we get the
following.
\proclaim{Corollary 3} If a norm $\|\cdot\|_X$ on $\r2$ is defined by a 
polygon $BX$ with finitely many sides all of which have algebraic slopes
then there is a compact $A\subset X$ such that the Hausdorff dimension of $A$
is $2$ and Lebesgue measure of $\Delta_X(A)$ is $0$. 
\endproclaim
It would be interesting to know if the result is true without supposition
on the slopes of the sides.

Recall that, by \cite{I{\L}03}, it is enough to prove the implication 
$(ii)\rightarrow(i)$ in Theorem 3 assuming that $BX$ is a polygon. In that case
we prove a stronger result.
\proclaim{Theorem 4}Let $BX$ be a polygon with finitely many sides which does
not satisfy the condition (i)
of Theorem 3. Then for any well-distributed set $S$ we have
$$|\Delta_{X,N}(S)|\gg N\log N/\log\log N\quad(N\ge3).\tag0.5$$
\endproclaim
Comparison of Theorem 4 with Theorem 2 shows that the growth of 
$|\Delta_{X,N}(S)|$ for well-distributed sets and $N\to\infty$ does not 
distinguish the spaces $X$ with small and big cardinality of $Ext(BX)$. 

\bigskip
\centerline{\S1. PROOF OF THEOREMS 1 AND 2}
\demo{Proof of (i)} Without loss of generality we may assume that 
$BX\subset Bl_2^2$ and the set $S$ is well-distributed in $X$ with the 
constant $K=1/2$. Also, choose $\delta>0$ so that
$$\delta Bl_2^2\subset BX.\tag1.1$$
By (0.1), for any $\ve>0$ there are arbitrary large $N_0$ such that
$$|\Delta_{X,N_0}(S)|\le\ve N_0^{3/2}.$$
If $N_0\ge8$ then the number of integers $j\ge0$ with $N_0/2+4j\le N_0-2$
is
$$\ge(N_0/2-2)/4\ge N_0/8.$$
Thus, there is at least one $j$ such that $N=N_0/2+4j$ satisfies the condition
$$|(\Delta_{X}(S)|\cap(N-2,N+2))\le8\ve N_0^{3/2}/N_0\le12\ve N^{1/2}.\tag1.2$$
So, (1.2) holds for arbitrary large $N$. 

We take any $N$ satisfying (1.2) and an arbitrary $P\in S$. Let $Q$ be
the closest point to $P$ in the space $X$ (observe that it exists since 
$S$ is closed due to (0.1)). Then, by well-distribution of $S$ (recall
that $K=1/2$) we have 
$$\|P-Q\|_X\le1.\tag1.3$$ 
Without loss of generality, $P=0$.
Denote $M=[2N\delta]$ and consider the rays
$$L_j=\{(r,\theta):\theta=\theta_j=2\pi j/M\},$$
where $(r,\theta)$ are the polar coordinates in $l_2^2$. Consider a 
point $R_j$, $1\le j\le M$, with the polar coordinates 
$(r_j,(\theta_{j-1}+\theta_j)/2)$ such that $\|R_j\|_X=N$. By (1.1) we have
$$r_j\ge\delta N.$$
Therefore, the Euclidean distance from $R_j$ to the rays $L_{j-1}$ and $L_j$
is
$$r_j\sin(\pi/M)\ge N\delta\sin(\pi/(2N\delta))>1.\tag1.4$$
provided that $N$ is large enough. Therefore, the distance from $R_j$ to these
rays in $X$ is also greater than $1$. Also, the distance from $R_j$ to the 
circles
$$\Gamma_1=\{R:\|R\|_X=N-1\},\quad \Gamma_2=\{R:\|R\|_X=N+1\}$$
in $X$ is equal to $1$. Thus, the $X$-disc of radius $1/2$ with the center at 
$R_j$ is contained in the open region $U_j$ bounded by $L_{j-1}$, $L_j$, 
$\Gamma_1$, and $\Gamma_2$. By the supposition on $S$ there is a point
$P_j\in U_j\cap S$.

Observe that for any $j$ we have 
$$N-1<\|P-P_j\|_X<N+1,\quad N-2<\|Q-P_j\|_X<N+2.$$
Let $U=\{(\|P-P_j\|_X,\|Q-P_j\|_X)\}$. By (1.2),
$$|U|\le144\ve^2N.\tag1.5$$
For any $(n_1,n_2)\in U$ we denote 
$$J_{n_1,n_2}=\{j:\|P-P_j\|_X=n_1,\quad \|Q-P_j\|_X=n_2\}.$$
By \cite{I{\L}03, Lemma 1.4, (i)}, if $j_1,j_2,j_3\in J_{n_1,n_2}$
then one of the points $P_{j_1},P_{j_2},P_{j_3}$ must lie on the segment
connecting two other points and contained in the circle $\{R:\|P-R\|_X=n_1\}$.
This implies that for all $j\in J_{n_1,n_2}$ but at most two indices
the intersection of $\p$ with the sector $S_j$ bounded by $L_{j-1}$ and $L_j$
is inside some line segment contained in $\p$. Therefore, by (1.5), the number of 
sectors $S_j$ containing an extremal point of $BX$ is at most $288\ve^2N$. 
For $R\in\p$ with the polar coordinates $(r,\theta)$ denote 
$\Theta(R)=\theta$. Define the measure on $\p$ in such a way that for any 
Borel set $V\subset \p$ the measure $\mu_P(V)$ is defined as the Lebesgue measure 
of $\Theta(V)$. In particular,
$$\mu_p(\p\cap S_j)=\frac{2\pi}M.$$
Clearly, $\mu_p$ is equivalent to the standard Lebesgue 
measure on $\p$. We have proved that
$$\mu_p(Ext(BX))\le288\ve^2N\frac{2\pi}M.$$
But $1/M\le1/(N\delta)$. Hence,
$$\mu_p(Ext(BX))\le2\pi\times288\ve^2/\delta.$$
As $\ve$ can be chosen arbitrarily small, we get $\mu_p(Ext(BX))=0$, and 
this completes the proof of (i).

Proof of (ii) follows the same scheme. Inequality (1.2) should be replaced by
$$|\Delta_{X}(S)\cap(N-2,N+2)|\le\Delta N^{\alpha},$$
where $\Delta$ may depend only on $X$, $S$, and $\alpha$. We define the
distance $d_p$ on $\p$ as the distance between the polar coordinates.
This metric is equivalent to the $X$-metric. The set $Ext(BX)$ can be covered
by at most $2\Delta^2N^{2\alpha}$ arcs $\p\cap S_j$ each of them has the 
$d_p$-diameter at most $2\pi/(N\delta)$. This implies the required estimate 
for the Hausdorff dimension of $Ext (BX)$. 
\enddemo  

\demo{Proof of Theorem 2} We select an increasing sequence $\{N_j\}$
of positive integers such that
$$\psi(N)\ge5^j\quad(N\ge N_j).\tag1.6$$
By $\Lambda_j$ we denote the set of numbers $a/q$ with $a\in\Bbb Z$, $q\in\n$,
$q\le N_j$. We will construct a ball $BX$ on the Euclidean plane. Moreover, 
it will be symmetric with respect to the lines $x_1=x_2$ and $x_1=-x_2$, and thus 
it suffices to construct $BX$ in the quadrant $Q=\{x:\,x_2\ge|x_1|\}$.

Let $D_0$ be the square
$$D_0=\{x:\,0\le x_2+x_1\le1,\,0\le x_2-x_1\le1\}.$$
We will construct a decreasing sequence of polygons $D_j$; each one
will be defined as a result of cutting some angles from the previous one.
The sides $V_1,V_2$ of $D_0$ with an endpoint at the origin will not be 
changed. The intersection of the sequence $D_j$ will define the part of our 
$BX$ in $Q$. In particular, the points $(\pm1/2,1/2)$ will be vertices of all
polygons $D_j$. Therefore, these points as well as the symmetrical points
$(\pm1/2,-1/2)$ will be in $\p$.

First, we construct $D_1$ as a result of cutting $D_0$ by a line $x_2=u$ for 
some $u\in(1/2,1)$. We choose
$u$ such that for intersection points $x^1$ and $x^2$ of this line with
the boundary of $D_0$ the ratios $x_1^j/x_2^j\,(j=1,2)$ differs from all 
numbers $\lambda\in\Lambda_1$. Moreover, we take 
neighborhoods $U_j$ of the points $x^j$ ($j=1,2$) such that
$$\forall y\in U_j\,\,y_2/y_1\not\in\Lambda_1\quad(j=1,2).$$
In the sequel we shall make other cuts only inside the sets $U_1$ and $U_2$.
This means that all points $x$ on the boundary of $D_1$ with 
$x_1/x_2\in\Lambda_1$ not belonging to the sides $V_1$, $V_2$ as well as their
neighborhoods in the boundary of $D_1$ will remain in all polygons 
$D_2,D_3,\dots$, and eventually they will be interior points of some segments 
in the boundary of $BX$ with a slope $-1$, $0$, or $1$, 

On the second step, we construct $D_2$ as a result of cutting $D_1$ by lines 
with slopes $-1/2$ and $1/2$ such that for any new vertex $x$ of a polygon 
$D_2$ we have $x_2/x_1\not\in \Lambda_2$. Moreover, we take 
neighborhoods $U(x)$ of all these points $x$ (each is contained in $U_1$
or in $U_2$) such that
$$\forall y\in U(x)\,\,y_2/y_1\not\in\Lambda_2.$$
Again, we shall make other cuts only inside the sets $U(x)$.
This means that all points $x$ on the boundary of $D_2$ with 
$x_1/x_2\in\Lambda_2$ not belonging to the sides $V_1$, $V_2$ as well as their
neighborhoods in the boundary of $D_2$ will remain in all polygons 
$D_3,D_4,\dots$, and eventually they will be interior points of some segments 
in the boundary of $BX$ with a slope $a/2$, $a\in\Bbb Z$, $|a|\le2$.

Proceeding in the same way, we shall get a ball $BX$ with the following 
property: if $x\in\p$ and $x_1/x_2\in \Lambda_{j+1}$ for some $j$ then
$x$ is an interior point of some segment contained in $\p$ with a slope
$a/2^j$, $a\in\Bbb Z$, $|a|\le2^j$. This segment is a part of
a line $2^jx_2-ax_1=b(a,j)$ or a symmetrical line 
$2^jx_2-ax_1=-b(a,j)$. Also, by symmetry, if $x\in\p$ and 
$x_2/x_1\in \Lambda_{j+1}$ for some $j$ then $2^jx_1-ax_2=b(a,j)$ or 
$2^jx_1-ax_2=-b(a,j)$. In terms of the norm $\|\cdot\|_X$ we conclude
that if $x\in X$ and $x_1/x_2\in\Lambda_{j+1}$ or $x_2/x_1\in\Lambda_{j+1}$ 
then $\|x\|_X$ is equal to one of the numbers $|2^jx_1-ax_2|/|b(a,j)|$ or 
$|2^jx_2-ax_1|/|b(a,j)|$, $a\in\Bbb Z$, $|a|\le2^j$. Also, observe that, by 
our construction, $BX$ is contained in the square $[-1,1]^2$. Therefore,
$$\|x\|_X\ge\max(|x_1|,|x_2|).\tag1.7$$ 

Now let us take the lattice $S=\z2$ and estimate $|\Delta_{X,N}(S)|$
for $N_j<N\le N_{j+1}$. If $x,y\in S$ and $\|x-y\|_X\le N$, then
we have $\|x-y\|_X=|(z_1,z_2)|_X$ where $z_1,z_2\in\Bbb Z$ and, by (1.7),
$\max(|z_1|,|z_2|)\le N$. Hence, $(z_1,z_2)=(0,0)$, or 
$x_1/x_2\in\Lambda_{j+1}$, or $x_2/x_1\in\Lambda_{j+1}$. Therefore,
$\|x-y\|_X=0$ or $\|x-y\|_X$ is equal to one of the numbers 
$|2^jx_1-ax_2|/|b(a,j)|$ or $|2^jx_2-ax_1|/|b(a,j)|$, $a\in\Bbb Z$, 
$|a|\le2^j$. For every $a$ we have 
$$|2^jx_1-ax_2|\le2^j|x_1|+|a|\times|x_2|\le2^{j+1}N.$$
Taking the sum over all $a$ we get
$$|\Delta_{X,N}(S)|\le(2^{j+1}+1)2^{j+1}N+1\le2^{2j+3}N.\tag1.8$$
On the other hand, by (1.6),
$$\psi(N)\ge5^j.\tag1.9$$
Comparing (1.8) and (1.9), we get (0.4) and thus complete the proof of the 
theorem.
\enddemo

\bigskip
\centerline{\S2. PROOF OF THEOREM 3, PART I} 
\bigskip

In this section we prove that the condition (i) of Theorem 3 implies (ii).

Assume that $\p$ consists of a finite number of line segments with slopes
$\bb_1,\bb_2,\dots,\bb_r$, all real and algebraic.  Let  
$\ff_\qq[\bb_1,\dots,\bb_r]$ be the field extension of $\qq$ 
generated by $\bb_1,\dots,\bb_r$, and let $\aa_0$ be its primitive
element, i.e. an algebraic number such that 
$\ff_\qq[\bb_1,\dots,\bb_r]=\ff_\qq[\aa_0]$.  
We may assume that $\aa_0$ is an algebraic integer: indeed, if $\aa_0$
is a root of $P(x)=a_dx^d+\dots+a_0$, then $\aa'_0=a_d\aa_0$ is a root
of $a_d^{d-1}P(x/a_d)=x^d+a_{d-1}x^{d-1}+a_{d-2}a_dx^{d-2}+\dots+a_0a_d^{d-1}$,
hence an algebraic integer, and generates the same extension.

It suffices to prove that
there is a well-distributed set $S\subset \r2$ such that 
$$
|\{x+\beta y:\ (x,y)\in S-S,\ |x|+|y|\leq R\}|\ll R, \tag2.1
$$
for each $\beta\in\ff_\qq[\alpha]$.

Since $\ff_\qq[\bb_1,\dots,\bb_r]\subset\rr$, we have $\aa_0\in\rr$.
Let $\aa_1,\dots,\aa_{d-1}$ be the algebraic conjugates of $\aa_0$ in ${\Bbb C}$
(of course they need not belong to $\ff_\qq[\aa_0]$).
Define for $C>0$
$$T(C)=\{\sum_{j=0}^{d-1} a_j\aa_0^j:\ a_j\in\zz,
|\sum_{j=0}^{d-1} a_j\aa_k^j|\leq C,\ k=1,\dots,d-1\},$$
and
$$S=T(C)\times T(C),$$
where $C$ will be fixed later.  

We first claim that $T(C)$ is well distributed in $\rr$ (with the implicit
constant dependent on $C$), and that 
$$|T(C)\cap[-R,R]|\ll R. \tag2.2$$  
Indeed, let $x=(x_0,x_1,\dots,x_{d-1})^T$ solve
$$
\sum_{j=0}^{d-1}\aa_0^j\,x_j=1,$$
$$\sum_{j=0}^{d-1}\aa_k^j\,x_j=0,\ k-1,\dots,d-1.$$
Since the Vandermonde matrix $A=(\aa_k^j)$ is nonsingular, ${x}$
is unique.  In particular, it follows that ${x}$ is real-valued;
this may be seen by taking complex conjugates and observing that 
$\aa_k$ is an algebraic conjugate of $\aa_0$ if and only if so is
$\bar\aa_k$, hence $\bar{x}$ solves the same system of equations.

To prove the first part of the claim, it suffices to show that 
there is a constant $K_1$ such that for any $y\in\rr$ there is a
$v\in T(C)$ with $|y-v|\leq K_1$.  Fix $y$, then we have 
$$y=\sum_{j=0}^{d-1}\aa_0^j\,yx_j.$$
Let $v_j$ be an integer such that $|v_j-yx_j|\leq 1/2$, and let
$v=\sum_{j=0}^d \aa_0^j v_j$.  Then
$$|y-v|=|\sum_{j=0}^{d-1}\aa_0^j(yx_j-v_j)|
\leq\frac{1}{2}\sum_{j=0}^{d-1}|\aa_0^j|=:K_1,$$
and, for $k=1,\dots,d-1$,
$$|\sum_{j=0}^{d-1}\aa_k^j\,v_j|
\leq |\sum_{j=0}^{d-1}\aa_k^j\,(yx_j-v_j)|
+y|\sum_{j=0}^{d-1}\aa_k^j\,x_j|
\leq \frac{1}{2}\sum_{j=0}^{d-1}|\aa_k^j|.$$
The claim follows if we let $C\geq \frac{1}{2}\sum_{j=0}^{d-1}|\aa_k^j|.$

We now prove (2.2).  It suffices to verify that
there is a constant $K_2$ such that for any $y\in\rr$ there are at most
$K_2$ elements of $T(C)$ in $[y-C,y+C]$.  Let 
$a=\sum_{j=0}^{d-1}\aa_0^j\,a_j$, then the conditions that $a\in T(C)$
and $|y-a|\leq C$ imply that
$$A\tilde a-\tilde y\in CQ,$$
where $\tilde a=(a_0,\dots,a_{d-1})^T$, $\tilde y=(y,0,\dots,0)^T$, 
and $Q=[-1,1]^d$.  In other words, $\tilde a\in A^{-1}\tilde y
+CA^{-1}Q$.  But it is clear that the number of integer lattice points contained
in any translate of $CA^{-1}Q$ is bounded by a constant.

It remains to prove (2.1).  
Observe first that if $x,x'\in T(C)$, then $x-x'\in T(2C)$.  Thus,
in view of (2.2), it is enough to prove
that for any two algebraic integers $\beta,\gamma\in\zz_\qq[\alpha]$ 
there is a $C_1=C_1(\beta,\gamma)$
such that if $x,y\in T(2C)$, then $x\beta+ y\gamma\in T(C_1)$.  By the
triangle inequality, it suffices to prove this with $y=0$.  Let
$x\in T(C)$, then $x=\sum_{j=0}^{d-1}\aa_0^j\,x_j$ for some $x_j\in
\zz$.  We also write $\beta=\sum_{j=0}^{d-1}\aa_0^j\,b_j$, with $b_j\in\zz$.
Then $\beta y=\sum_{i,j=0}^{d-1}\aa_0^{i+j}\,x_ib_j$.
We thus need to verify that
$$|\sum_{i,j=0}^{d-1}\aa_k^{i+j}\,x_ib_j|\leq C_1$$
for $k=1,\dots,d-1$.  But the left side is equal to
$$|\sum_{i=0}^{d-1}\aa_k^{i}\,x_i|
\cdot |\sum_{j=0}^{d-1}\aa_k^{j}\,b_j|,$$
which is bounded by $C_1(\beta)=C
\max_k |\sum_{j=0}^{d-1}\aa_k^{j}\,b_j|$.

\bigskip

{\bf Example.} Let $BX$ be a symmetric convex octagon whose 
sides have slopes $0,-1,\infty, \sqrt{2}$.  Let also $T(C)
=\{i+j\sqrt{2}:\ |i-j\sqrt{2}|\leq C\}$, and $S=T(10)\times T(10)$.
It is easy to see that $T(C)$ is well distributed and that 
(2.2) holds.
Let $x,y\in S$, then $x-y=(i+j\sqrt{2},k+l\sqrt{2})$, where 
$i+j\sqrt{2},k+l\sqrt{2}\in T(20)$.  Depending on where 
$x-y$ is located, the distance from
$x$ to $y$ will be one of the following numbers:
$$c_1|i+j\sqrt{2}|,$$
$$c_2|k+l\sqrt{2}|,$$
$$c_3|(i+k)+(j+l)\sqrt{2}|,$$
$$c_4|(i+j\sqrt{2})\sqrt{2}-(k+l\sqrt{2})|
=c_4|(2j-k)+(i-l)\sqrt{2}|.$$
Clearly, the first three belong to $T(20\max(c_1,c_2,c_3))$.
For the fourth one, we have
$$
c_4|(2j-k)-(i-l)\sqrt{2}|
=c_4|-(i-j\sqrt{2})\sqrt{2}-(k-l\sqrt{2})|$$
$$\leq 20 c_4(1+\sqrt{2}).$$
Hence all distances between points in $S$ belong to $T(C)$ for
some $C$ large enough, and in particular satisfy the cardinality
estimate (2.2).

\bigskip

\centerline{\S3. ADDITIVE PROPERTIES OF MULTIDIMENSIONAL SETS}
\centerline{AND SETS WITH SPECIFIC ADDITIVE RESTRICTIONS}

\bigskip
Let $Y$ be a linear space over $\Bbb R$ or over $\q$. For $A,B\subset Y$ and 
$\aa\in\rr$ or $\qq$ we denote
$$A+B=\{a+b:\,a\in A,b\in B\},\ \aa A=\{\alpha a:\,a\in A\} .$$
We say that a set $A\subset Y$ is a $d$-dimensional if $A$ is contained in 
some $d$-dimensional affine subspace of $Y$, but in no $d-1$-dimensional 
affine subspace of $Y$.  We will denote the dimension of a set $A$ 
by $d_A$.

The following result is due to Ruzsa \cite{Ru94, Corollary 1.1}.

\proclaim{Lemma 3.1} Let $A,B\subset\rr^d$, $|A|\leq|B|$, and assume
that $A+B$ is $d$-dimensional.  Then
$$|A+B|\ge |B|+d|A|-d(d+1)/2.\tag3.1$$
\endproclaim

The special case of Lemma 3.1 with $A=B$ was proved earlier
by Freiman \cite{F73, p.~24}). In this case we also have the
following corollary.

\proclaim{Corollary 3.1} Let $A\subset\rr^d$, and assume that
$|A+A|\leq K|A|$, $K\leq |A|^{1/2}$.  Then the dimension of $A$ 
does not exceed $K$.
\endproclaim

\demo{Proof}
Let $|A|=N\geq 1$, then $d_A\leq N-1$. Suppose that $d_A> K$.  The function 
$f(x)=(x+1)N-x(x+1)/2$ is increasing for $x\leq N-1/2$, hence by (3.1)
we have
$$KN\geq f(d_A)> f(K)=(K+1)N-\frac{K(K+1)}{2},$$
i.e. $K(K+1)>2 N$, which is not possible if $K^2\leq N$.
\enddemo

\bigskip

We observe that Lemma 3.1, and hence also Corollary 3.1,
 extends to the case when $A,B$ are subsets of
a linear space $Y$ over $\q$. Assume that $Y$ is $d$-dimensional, and
take a basis $\{e_1,\dots.e_d\}$ in $Y$. 
Consider the space $\Bbb R^d$ with a basis $\{e_1',\dots.e_d'\}$. 
We can arrange a mapping $\Phi:Y\to Y'$ by
$$\Phi(\sum_{j=1}^d\alpha_j e_j)=\sum_{j=1}^d\alpha_j e_j'.$$
It is easy to see that $\Phi$ is Freiman's isomorphism of any order
and, in particular, of order $2$: this means that for any
$y_1,,y_2,z_1,z_2$ from $Y$ the condition
$$y_1+y_2\neq z_1+z_2$$
implies
$$\Phi(y_1)+\Phi(y_2)\neq \Phi(z_1)+\Phi(z_2).$$
Therefore, if $A,B$ are finite subsets of $Y$ and $A'=\Phi(A), B'=\Phi(B)$,
then $|A+B|=|A'+B'|$,
and we get the required inequality for $|A+B|$. 

The following is a special case of \cite{N96, Theorem 7.8}.

\proclaim{Lemma 3.2} If $N\in\n$, $K>1$, $A\subset Y$, and $B\subset Y$ 
satisfy
$$\min(|A|,|B|)\ge N,\quad |A+B|\le KN,\tag3.2$$
we have
$$|A+A|\le K^2|A|.$$
\endproclaim

\proclaim{Corollary 3.2}
If $N\in\n$, $K>1$, and if $A,B\subset Y$ satisfy (3.2) for some
K with $K^2(2K^2+1)<N$, then $d_{A+B}\leq K$.  In particular, 
$d_A\leq K$ and $d_B\leq K$.
\endproclaim

\demo{Proof}
By Lemma 3.2, we have $|A+A|\leq K^2N$, hence Corollary 3.1 implies
that 
$$d_A\leq K^2,$$
and similarly for $B$.  Hence $d_{A+B}\leq d_A+d_{B}\leq 2K^2$.
By Lemma 3.1, we have
$$KN\geq |A+B|\geq (1+d_{A+B})N-\frac{d_{A+B}(d_{A+B}+1)}{2}$$
$$\geq d_{A+B}N+N-K^2(2K^2+1)\geq d_{A+B}N,$$
which proves the first inequality.  To complete the proof, 
observe that $d_{A+B}\geq \max(d_A,d_B)$.
\enddemo

\proclaim{Lemma 3.3} Let $K>0$, $A$ and $B$ be finite
nonempty subsets of $\Bbb R$, $\alpha\in\Bbb R\setminus\{0\}$. Also, suppose 
that the following conditions are satisfied
$$|A-\alpha B|\le K|B|.\tag3.3$$
Then there is a set $B'\subset B$ such that
$$|A-\alpha B'|\le K|B'|,\tag3.4$$
$$|B'|\ge|A|/K,\tag3.5$$
and for any $b_1,b_2\in B'$ the number $\alpha(b_1-b_2)$ is a linear 
combination of differences $a_1-a_2$, $a_1,a_2\in A$, with integer 
coefficients.
\endproclaim
\demo{Proof} Let us construct a graph $H$ on $B$. We join $b_1,b_2\in B$ 
(not necessary distinct) by an edge if there are $a_1,a_2\in A$ such that 
$a_1-\alpha b_1=a_2-\alpha b_2$. Let $B_1,\dots,B_s$ be the components of 
connectedness of the graph $H$. Thus, for any $j=1,\dots,s$ and for any 
$b_1,b_2\in B_s$ there is a path connecting $b_1$ and $b_2$ and consisting
of edges of $H$ (a one-point path for $b_1=b_2$ is allowed). This implies
that $\alpha(b_1-b_2)$ is a sum of differences $a_1-a_2$ for some pairs 
$(a_1,a_2)\in A\times A$. Also, denoting
$$S=A-\alpha B,\quad S_j=A-\alpha B_j,$$
we see that, by the choice of $B_1,\dots,B_s$, the sets $S_j\,(j=1,\dots,s)$ 
are disjoint. 

Since
$$|B|=\sum_{j=1}^s|B_j|,\quad|S|=\sum_{j=1}^s|S_j|,$$
there is some $j$ such that
$$|S_j|/|B_j|\le|S|/|B|,$$
and, by (3.3),
$$|S_j|\le K|B_j|.$$ 
On the other hand,
$$|S_j|=|A-\alpha B_j|\ge|A|.$$
Hence,
$$|B_j|\ge|S_j|/K\ge|A|/K.$$
So, the set $B'=B_j$ satisfies (3.4) and (3.5), and Lemma 3.3 follows.
\enddemo

\proclaim{Lemma 3.4} Let $K>0$, $A$ and $B$ be finite
nonempty subsets of $\Bbb R$, $\alpha_1,\alpha_2\in\Bbb R\setminus\{0\}$. 
Also, suppose that the conditions 
$$|A-\alpha_1B|\le K|B|,\tag3.6$$
$$|A-\alpha_2B|\le K|A|,\tag3.7$$
are satisfied. Then there are nonempty sets $A'\subset A$ and 
$B'\subset B$ such that
$$|A-\alpha_1 B'|\le K|B'|,\tag3.8$$
$$|A'-\alpha_2 B'|\le K|A'|,\tag3.9$$
$$|A'|\ge|A|/K^2,\tag3.10$$
and for any $a_1',a_2'\in A'$ the difference $a_1'-a_2'$ is a linear 
combination of numbers $\frac{\alpha_2}{\alpha_1}(a_1-a_2)$, $a_1,a_2\in A$, 
with integer coefficients.
\endproclaim
\demo{Proof} By (3.6), we can use Lemma 3.3 for $\alpha=\alpha_1$, and we get
(3.8) and (3.5). Further, we use Lemma 3.3 again for $B',A$ (thus, in the
reverse order), and we get (3.9) and also
$$|A'|\ge|B'|/K.$$
Combining the last inequality with (3.5) we obtain (3.10). 
The proof of the lemma is complete.
\enddemo

Replacing (3.8) by a weaker inequality
$$|A'-\alpha_1 B'|\le K|B'|$$
and iterating Lemma 3.4, we get the following.
\proclaim{Lemma 3.5} Let $K>0$, $A$ and $B$ be finite
nonempty subsets of $\Bbb R$, $\alpha_1,\alpha_2\in\Bbb R\setminus\{0\}$. 
Also, suppose that the conditions (3.6) and (3.7) 
are satisfied. Then there are nonempty sets $A_j\subset A$ and $B_j\subset B$ 
($j=0,1,\dots,$) such that $A_0=A$, $B_0=B$, $A_j\subset A_{j-1}$,
$B_j\subset B_{j-1}$ for $j\ge1$,
$$|A_j-\alpha_2B_j|\le K|A_j|\quad(j\ge1),$$
$$|A_j|\ge|A|/K^{2j},$$
and for any $a_1,a_2\in A_j$ the difference $a_1-a_2$ is a linear 
combination of numbers $\frac{\alpha_2^j}{\alpha_1^j}(a_1'-a_2')$, 
$a_1',a_2'\in A$, with integer coefficients.
\endproclaim
Now we are in position to come to the main object of our constructions:
to show that under the assumptions of Lemma 3.5, providing that the number
$\alpha_1/\alpha_2$ is transcendental, we can conclude that the dimension
of the set $A$ over $\q$ cannot be too small. 
\proclaim{Corollary 3.6} Let $K>0$, $A$ and $B$ be finite nonempty 
subsets of $\Bbb R$, $\alpha_1,\alpha_2\in\Bbb R\setminus\{0\}$ 
such that $\alpha_1/\alpha_2$ is transcendental. Also, suppose that the 
conditions (3.6) and (3.7) are satisfied. Then, if for some $d\in\n$ the 
inequality
$$|A|>K^{2d}\tag3.11$$
holds, then the dimension of $A$ over $\q$ is greater than $d$.
\endproclaim
\demo{Proof} By Lemma 3.5 and (3.11), we have $|A_d|\ge2$. Take distinct 
$a_1,a_2\in A_d$. Then also $a_1,a_2\in A_j$ for $j=0,1,\dots,d$, and, by 
Lemma 3.6, the difference $a_1-a_2$ is a linear combination of numbers 
$\frac{\alpha_2^j}{\alpha_1^j}(a_1'-a_2')$, $a_1',a_2'\in A$, with integer 
coefficients. Therefore, all numbers 
$b_j=\frac{\alpha_1^j}{\alpha_2^j}(a_1-a_2)$ belong
to the linear span of $a_1'-a_2'$, $a_1',a_2'\in A$, over $\q$. But, since
$\alpha_1/\alpha_2$ is transcendental, the numbers $b_j\,(j=0,\dots,d)$
are linearly independent over $\q$. Therefore, the dimension of the linear 
span of $a_1'-a_2'$, $a_1',a_2'\in A$, over $\q$ is at least $d+1$, as
required.
\enddemo

\proclaim{Corollary 3.7} If $A$ is a subset of $\Bbb R$, $2\le|A|<\infty$,
$\alpha$ is a transcendental real number, then
$$|A-\alpha A|\gg|A|\log|A|/\log\log|A|.$$ 
\endproclaim

\demo{Proof} Suppose that the conclusion fails, then 
for any $\eps>0$ we may find arbitrarily large $N$ and $A\subset
\rr$ with $|A|=N$ such that 
$$|A-\aa A|\leq KN,\ K=\eps\frac{\log N}{\log\log N}.$$
By Corollary 3.2, we have $d_A\leq K$.  On the other hand,
(3.6) holds with $B=A$, $\alpha_1=\alpha$, and, since $A-\alpha^{-1}A
=-\alpha^{-1}(A-\alpha A)$, (3.7) holds with $B=A$ and $\aa_2=
\alpha^{-1}$.  
Corollary 3.7 then implies that
$$N\le K^{2K}.$$
Taking logarithms of both sides, and assuming that $2\eps<1$, we obtain
$$\log N\leq 2\eps\frac{\log N}{\log\log N}
(\log(2\eps)+\log\log N-\log\log\log N)\leq 2\eps\log N,$$
which is not possible if $N$ was chosen large enough.
\enddemo

{\bf Remark.} On the other hand, if $\alpha\in\Bbb R$ is an algebraic number,
then one can use our construction from \S2 to show that for any $N\in\n$
there is a set $A\subset\Bbb R$, $|A|=N$, such that
$$|A-\alpha A|\le C|A|,$$
where $C$ depends only on $\alpha$.

Finally, we state a lemma due to J.~Bourgain\cite{B99, Lemma 2.1}.  For our
purposes, we need a slightly more precise formulation than that given in
\cite{B99}; the required modifications are described below.

\proclaim{Lemma 3.8} Let $N\ge2$, $A,B$ be finite subsets of $\rr$ and 
$G\subset A\times B$ such that
$$|A|,|B|\le N,\tag3.12$$
$$|S|\le N\quad\text{where}\quad S=\{a+b:(a,b)\in G\},\tag3.13$$
$$|G|\ge\delta N^2.\tag3.14$$
Then there exist $A'\subset A$, $B'\subset B$ satisfying the conditions
$$|(A'\times B')\cap G|\gg\delta^5N^2(\log N)^{-C_1},\tag3.15$$
$$|A'-B'|\ll N^{-1}(\log N)^{C_2}\delta^{-13}|(A'\times B')\cap G|.\tag3.16$$
\endproclaim

In \cite{B}, the bounds (3.15) and (3.16) involved factors
of the form $N^{\gamma+}$ and $N^{\gamma-}$, where $N^{\gamma+}$ ($N^{\gamma-}$)
means $\le C(\ve)N^{\gamma+\ve}$ for all $\ve>0$ and some $C(\ve)>0$ (resp., 
$\ge c(\ve)N^{\gamma-\ve}$ for all $\ve>0$, $c(\ve)>0$). We need 
a slightly stronger statement, namely that the same bounds hold with the 
factors in question obeying the
inequalities $\ll N^{\gamma}(\log N)^C$ or $\gg N^{\gamma}(\log N)^{-C}$,
respectively, for some appropriate choice of a constant $C$. 
A careful examination of the proof in \cite{B99} shows that it remains valid 
with this new meaning of the notation $N^{\gamma+}$ and $N^{\gamma-}$,
and that one may in fact take $C_1=5$, $C_2=10$.
We further note that although Bourgain states his lemma for $A,B\subset
\zz^d$, the same proof works for $A,B\subset\rr$ if the exponential sum
inequality \cite{B99,(2.7)} is replaced by
$$|{G}|<\int_S{\chi_A*\chi_B}\leq|S|^{1/2}\|\chi_A*\chi_B\|_2;$$
we then observe that 
$$\|\chi_A*\chi_B\|_2^2 
=|\{(a,a',b,b')\in A\times A\times B\times B: a+b=a'+b'\}|$$
$$=|\{(a,a',b,b')\in A\times A\times B\times B: a-b'=a'-b\}|
=\|\chi_A*\chi_{-B}\|_2^2,$$
and proceed further as in \cite{B99}.  A similar modification should
be made in\newline
\cite{B99,(2.36)}.

\bigskip

\centerline{\S4. PROOF OF THEOREM 4}

\bigskip

In this section we prove Theorem 4; note that this also proves the
implication (ii)$\Rightarrow$(i) of Theorem 3.

Suppose that $BX$ is a polygon with finitely many sides for which
the conclusion of the theorem fails, i.e. that there
is a well distributed set $S$ such that for any $\eps>0$ there
is an increasing sequence of positive integers $N_1,N_2,\dots\to
\infty$ with
$$|\Delta_{X,N_j}(S)|<\eps N_j\psi(N_j),\tag4.1$$
where
$$\psi(N)=\log N/\log\log N.$$
Without loss of generality we may assume that $\partial BX$ contains
a vertical line segment and a horizontal line segment, and that
$c_1Bl_2^2\subset BX\subset Bl_2^2$.  Let also $c_2\in(0,1/10)$ be a small
constant such that all sides of $BX$ have length at least $8 c_2$.

Let $M$ be a sufficiently large number which may depend on $\eps$; all
other constants in the proof will be independent of $\eps$. 
Let $T=N_{j_0}$ for some $j_0$ large enough so that $T>M$, and 
let $N=c_2T$.  Suppose that one of the two vertical sides of $BX$
is the line segment $\{(x_1,x_2):\ x_1=v_1,|x_2-v_2|\leq r\}$, where $v_1>0$.
Let also $Q=Int\,(N\cdot BX)$, $v=(v_1,v_2)$, and 
$$A=\{x_1:\ (x_1,x_2)\in S\cap Q\hbox{ for some }x_2\},$$
$$Q'=Q+(T-2N)v.$$
Observe that both $Q$ and $Q'$ have Euclidean diameter $\leq 2N$,
and that
$$Q'\subset\{(x_1,x_2): (T-3N)v_1 < x_1 < (T-N)v_1\},$$
so that 
$$\|x-x'\|_X\geq (1-4c_2)T> T/2, \ x\in Q,x'\in Q'.$$
By our choice of $c_2$ we have $c_2\leq r/4$, so that
$$T/2\cdot r\geq 2N.$$
Hence all $X$-distances between points in $Q$ and $Q'$ are measured
using the vertical segments of $\partial BX$, i.e.
$$\|x-x'\|_X=|x_1-x'_1|/v_1, \ x=(x_1,x_2)\in Q,x'=(x'_1,x'_2)\in Q_t.$$

Next, we claim that 
$$|\{\|x-x'\|_X:\ x\in S\cap Q,x'\in S\cap Q'\}|< K_0\eps N\psi(N),\tag4.2$$
where $K_0$ is a constant depending only on $c_2$.
Indeed, we have
$$\{\|x-x'\|_X:\ x\in Q,x'\in Q'\}
\subset [0,T],$$
hence the failure of (4.2) would imply that 
$$|\Delta_{X,T}(S)|\geq K_0\eps N\psi(N)
\geq \eps T\psi(T),$$
if $K_0$ is large enough (at the last step we used that $\psi(N)\gg
\psi(c_2^{-1}N)=\psi(T)$).  But this contradicts (4.1).

It follows that if we define
$$A'=\{x'_1:\ (x'_1,x'_2)\in S\cap Q'\hbox{ for some }x'_2\},$$
then we can estimate the cardinality of the difference set $A-A'$ using
(4.2):
$$|A-A'|<K_0\eps N\psi(N).\tag4.3$$
On the other hand, since $S$ is well distributed, we must have
$$|A|,|A'|\gg N.\tag4.4$$
Hence by Corollary 3.2 we have
$$d_A\ll \eps\psi(N).\tag4.5$$

We may now repeat the same argument with the vertical side of $\partial BX$
replaced by its other sides.  In particular, using the horizontal 
segment in $\partial BX$ instead, we obtain the following.  Let
$$B=\{x_2:\ (x_1,x_2)\in S\cap Q\hbox{ for some }x_1\},$$
then there is a set $B'\subset\rr$ such that
$$|B|,|B'|\gg N,\tag4.6$$
$$|B-B'|<K_0\eps N\psi(N),\tag4.7$$
$$d_B\ll \eps\psi(N).\tag4.8$$
Furthermore, assume that $\partial BX$ contains a segment of a line
$x_1+\alpha x_2=\beta$, then 
$$|\{x_1+\alpha x_2:\ (x_1,x_2)\in S\cap Q\}|\leq K_0\eps N\psi(N);\tag4.9$$
this estimate is an easier analogue of (4.3) obtained by counting
distances between points in $Q$ and just one point in the appropriate
analogue of $Q'$.

Suppose that $\partial BX$ contains segments of lines $x_1+\alpha_1x_2=C_1$,
$x_2+\alpha_2x_2=C_2$ (i.e. with slopes $-1/\alpha_1$, $-1/\aa_2$), where
$\aa_1,\aa_2$ are neither 0 nor $\infty$, 
and that the ratio $\alpha_1/\alpha_2$ is transcendental.  
Let $G=(A\times B)\cap S$,
then $|G|\geq c_4 N^2$ since $S$ is well distributed.  By (4.4), (4.6), and
(4.9) with $\aa=\aa_1$,
the assumptions of Lemma 3.8 are satisfied with $N$ replaced by 
$K_0\eps N\psi(N)$ and $\delta=c_4(K_0\eps \psi(N))^{-2}$.  We conclude that there
are subsets $A_1\subset A$ and $B_1\subset B$ such that
$$|(A_1\times B_1)\cap G|\gg N^2 \eps^c(\log N)^{-c},\tag 4.10$$
$$|A_1-\aa_1B_1|\ll N^{-1}\eps^{-c}(\log N)^{c}|(A_1\times B_1)\cap G|.\tag4.11$$
Here and below, $c$ denotes a constant which may change from line to line 
but is always independent of $N$.  We also simplified the right sides
of (4.10) and (4.11) by noting that $\psi(N)\leq \log N$.

Similarly, applying Lemma 3.8 with $G$ replaced by  $(A_1\times B_1)\cap G$
and $\aa_1$ replaced by $\aa_2$, we find
subsets $A_2\subset A_1$ and $B_2\subset B_1$ such that
$$|(A_2\times B_2)\cap G|\gg N^2\eps^{c}(\log N)^{-c},\tag 4.12$$
$$|A_2-\aa_2B_2|\ll N^{-1}\eps^{-c}(\log N)^{c}|(A_2\times B_2)\cap G|.\tag4.13$$
Clearly, (4.11) also holds with $A_1,B_1$ replaced by $A_2,B_2$.

Thus $A_2,B_2$ satisfy the assumptions (3.14), (3.15) of Corollary 3.7,
with $K=\eps^{-c}(\log N)^c$.  By (4.4), (4.5) and Corollary 3.7, we must have
for some constants $c,K_2$, 
$$cN\leq |A_2|<(\eps^{-1}\log N)^{K_2\eps\log N/\log\log N},$$
hence
$$\log c+\log N\leq \frac{K_2\eps\log N}{\log\log N}(\log\log N-\log\eps)
\leq 2K_2\eps\log N,$$
a contradiction if $\eps$ was chosen small enough.
This proves that if (0.5) fails, then the ratio between any two slopes,
other than $0$ or $\infty$, of sides of $BX$ is algebraic.

To conclude the proof of the theorem, we first observe that if $BX$ is
a rectangle, then there is nothing to prove.  If $BX$ is a hexagon
with slopes $0,\infty,\alpha$, we may always find a coordinate system
as in Theorem 3 (i); namely, if we let
$$x'_1=x_1,\ x'_2=\alpha x_2,\tag4.14$$
then the slopes $0$ and $\infty$ remain unchanged, and lines
$\alpha x_1- x_2=C$ with slope $\alpha$ are mapped to lines $x'_1-x'_2
=C/\alpha$ with slope $1$.  Finally, suppose that $BX$ is a polygon
with slopes $0,\infty,\alpha_1,\alpha_2,\dots,\alpha_l$, and apply
the linear transformation (4.14) with $\alpha=\alpha_1$.  Then
the sides of $\partial BX$ with slope $\alpha_1$ is mapped to 
line segments with slope 1; moreover, since the ratios $\alpha_j/\alpha_1$,
$j=2,3,\dots,l$, remain unchanged in the new coordinates, and since
we have proved that these ratios are algebraic, all remaining sides
of $\partial BX$ are mapped to line segments with algebraic slopes.

\bigskip

{\bf Acknowledgements.} This work was completed while the first author
was a PIMS Distinguished Chair at the University of British Columbia,
and was partially supported by NSERC grant 22R80520.
We are indebted to Ben Green for pointing out to us the reference [Ru].

\bigskip
\centerline{REFERENCES}
\bigskip

[B99] J.~Bourgain, On the dimension of Kakeya sets and related maximal inequalities,
Geom. Funct. Anal. {\bf 9} (1999), 256--282.

[E46] P. Erd\H{o}s, On sets of distances of $n$ points, Amer. Math. Monthly
{\bf 53} (1946), 248--250.

[F73] G. Freiman, Foundations of a structural theory of set addition (translation from
Russian), Translations of Mathematical Monographs, vol. 37, American Mathematical
Society, Providence, RI, 1973.

[I01] A. Iosevich, Curvature, combinatorics and the Fourier transform, Notices
Amer. Math. Soc. {\bf 46} (2001), 577--583.

[I{\L}03] A.~Iosevich and I.~{\L}aba, Distance sets of well-distributed planar 
point sets, Discrete Comput.~Geometry {\bf 31} (2004), 243--250.

[I{\L}04] A.~Iosevich and I.~{\L}aba, $K$-distance sets, Falconer conjecture and
discrete analogs, preprint, 2003.

[KT04] N.H.Katz and G.Tardos, A new entropy inequality for the Erd\H os
distance problem. in: Towards a Theory of Geometric Graphs.(ed.J Pach)
Contemporary Mathematics, vol. 342, Amer.Math Soc. 2004

[N96] M.~Nathanson, Additive Number Theory, II: Inverse Problems and the
Geometry of Sumsets, Springer-Verlag, New York, 1996.

[Ru94] I.~Ruzsa, Sum of sets in several dimensions, Combinatorica {\bf 14}
(1994), 485--490.

[ST01] J. Solymosi and Cs. T\'oth, Distinct distances in the plane, Discrete
Comput.~Geometry {\bf 25} (2001), 629--634.

[W99] T.~Wolff, Decay of circular means of Fourier transforms of measures,
Int. Math. Res. Notices {\bf 10} (1999), 547--567.

\bigskip

\enddocument